\documentclass[10pt]{article}

\usepackage{array}
\usepackage{amsmath}
\usepackage{amsfonts}
\usepackage{amssymb}
\usepackage{amsxtra}
\usepackage{xcolor}
\usepackage{epsfig}
\usepackage{multicol}
\usepackage{enumerate}

\title{Mean Values at Hopf Points and Oscillation-Induced Gain Modulation}

\author{William Harold Nesse \thanks{Department of Mathematics, University of Utah, Salt Lake City, UT (nesse@math.utah.edu).   }
\and Cooper John Hutchinson  \thanks{University of Utah, Salt Lake City, UT
   (Cooper.hutchinson@utah.edu).}   
}

\begin{document}

\maketitle

\begin{abstract}
We present a result concerning the mean value of orbits emerging from Hopf bifurcations. We then apply this result to identify a new phenomenon termed {\it oscillation-induced gain modulation}. A codimension one Hopf bifurcation of a system $\dot{x} = f(x; \alpha)$ with parameter $\alpha$ is characterized by the emergence of a limit cycle with an amplitude increasing from zero, coinciding with a change in the stability of an equilibrium $x_0(\alpha)$ when $\alpha$ passes through a critical value $\alpha^*$. This bifurcation is associated with the real part single eigenpair $\lambda = \mu(\alpha) \pm i \omega(\alpha)$ of the linearized system crossing zero: $\mu(\alpha^*) = 0$, $\mu'(\alpha^*) \neq 0$. We establish a result concerning the temporal mean of the oscillation cycle over the period $T$ of oscillation:  $\langle x \rangle_{\alpha} = \frac{1}{T} \int_0^{T} x(t; \alpha) dt $. We set the mean to be $\langle x \rangle_{\alpha}  = x_0(\alpha)$ when the equilibrium has no surrounding limit cycle. However, when a limit cycle exists, we show that that the deviation of the mean from the equilibrium is expressible as  $ \langle x \rangle_{\alpha}  - x_0(\alpha) = K \mu(\alpha) + \mathcal{O}(\mu(\alpha)^2)$. That is, the mean value deviates from the equilibrium's location in proportion to $\mu(\alpha)$, with a mean deviation determined by the vector quantity $K(\alpha) \mu(\alpha) $ that depends on the tensors of $f$ up to third-order. This result complements the well-known Hopf bifurcation property that the oscillation amplitude scales as  $\sqrt{\mu(\alpha)}$. If we consider $\alpha$ to be an input to the model, and the mean $\langle x \rangle_{\alpha} $ as the output, then the mean deviation $K \mu(\alpha)$  introduces a discontinuity to the cycle mean gain $\frac{d \langle x \rangle_{\alpha}}{d\alpha}$ at the bifurcation, which we term oscillation-induced gain modulation (OIGM). We prove this cycle mean deviation result for general Hopf points in two-dimensional and $n$-dimensional systems, as well as showcase several examples of OIGM.  
\end{abstract}



\section{Introduction}


Hopf bifurcations are an important mechanism to induce limit cycle solutions in parameterized dynamical systems \cite{hopf1942}. Hopf instabilities have widespread applications in the sciences, including in models of turbulent fluid flow \cite{hopf1948turbulence, caton1999, duvsek1994}; in the neurosciences, including in the Hodgkin-Huxley model \cite{hodgkin1952a,hodgkin1952b,hodgkin1952c,hodgkin1952d} \cite{jack1975electric}, as well as neural mass models \cite{wilson1972, powanwe21}; also in models of chemical reactions,  such as the Brusselator that characterizes the Belousov-Zhabotinsky reaction \cite{belousov1958, prigogine1968}; and additionally, models of multi-species population dynamics \cite{holling1965, may1972, zeeman1993hopf}, to name a few.
 
A simple Hopf bifurcation is defined as follows: Consider $n$-dimensional ($n\geq 2$) dynamical system $\dot{x} = f(x; \alpha)$ with parameter $\alpha$. Let $x_0(\alpha)$ be an equilibrium point ($f(x_0(\alpha); \alpha) = 0$) that exists for $\alpha$-values in some open interval. A Hopf instability is defined by a single pair of complex eigenvalues $\lambda$, $\bar{\lambda}$ of the linearized system, expanded about an equilibrium point $x_0(\alpha)$, crossing to the positive real half of the complex plane at a critical parameter value $\alpha^*$:  $\lambda(\alpha) = \mu(\alpha) \pm i \omega(\alpha)$, with $\mu(\alpha^*) = 0$. We also require that $\omega(\alpha) >0$, as well as a transversality condition $\mu'(\alpha^*) \neq 0$ to ensure non-degeneracy. When $\alpha$ passes through the $\alpha^*$-boundary into the interval $\tilde{I}_{\alpha^*}$, a single limit cycle solution $x(t)$ is assured to emerge surrounding the equilibrium $x_0(\alpha)$ with amplitude increasing from zero for $\alpha \geq \alpha^*$, and with period near $T \sim 2\pi/\omega$ ($x(t+T) = x(t)$), provided certain conditions on the third- and lower-order tensors of $f$ hold (i.e., a non-zero Lyapunov coefficient) \cite{hopf1942, marsden2012hopf}. 

This letter obtains a new result concerning the mean value of the system's limit-cycle solutions.  We define the mean value vector $\langle x \rangle_{\alpha}$ to be the temporal average over the limit cycle when an oscillation is present, but equal to the equilibrium $x_0(\alpha)$ otherwise:
\begin{align}
\label{xmean0}
\langle x \rangle_{\alpha} \equiv \begin{cases}  \frac{1}{T} \int_0^T x(t) dt, &  \alpha \in \tilde{I}_{\alpha^*} \\ x_0(\alpha),  &  \alpha \in I_{\alpha^*}  \end{cases},
\end{align}
where $ I_{\alpha^*}$ is the compliment of $ \tilde{I}_{\alpha^*}$---assume that $\alpha^*$ is included in $\tilde{I}_{\alpha^*}$ at its closed boundary. In this letter, we will focus on examples with {\it supercritical} Hopf bifurcations in which loss of equilibrium stability coincides with emergence of a stable limit cycle orbit; however, our results generalize to subcritical Hopf points as well. 


The main objective of this paper is to determine the deviation of the mean $\langle x \rangle_{\alpha}$ from the equilibria  $x_0(\alpha)$ when $\alpha$ crosses the bifurcation threshold. We will establish that the deviation between the mean and the equilibrium $\langle x \rangle_{\alpha} - x_0(\alpha)$ can be expressed by a $\mathcal{O}( \mu(\alpha))$ scaling:
 \begin{align}
 \label{mainResult0}
 \langle x \rangle_{\alpha} -  x_0(\alpha)  =  \begin{cases}   K(\alpha)\mu(\alpha), &  \alpha \in \tilde{I}_{\alpha^*} \\ 0,  &  \alpha \in I_{\alpha^*}  \end{cases},
\end{align}
where the vector quantity $K(\alpha)$ depends on, to lowest order, the first- through third-order tensors of $f$, and determines the relative direction the mean (\ref{xmean0}) to the equilibrium as a function of $\alpha$.  We will show examples where this vector function $K(\alpha)$ in (\ref{mainResult0}) can be zero, or non-zero, depending on properties of the model in question, and can produce dramatic, even orthogonal, trajectories relative to that of the equilibrium $x_0(\alpha)$.

This result (\ref{mainResult0}) has significance in many applications, owing to the ubiquity of Hopf-induced oscillatory instabilities found in models of natural systems. In particular, if $\alpha$ is considered to be an input to the system, and the gain of the mean response $ \frac{d\langle x \rangle_{\alpha}}{d\alpha}$ is an output of interest, this result (\ref{mainResult0}) shows that the  bifurcation will induce a discontinuity in the gain and consequently as sudden shift in the trajectory of the mean $ \langle x \rangle_{\alpha} $. We term this newly identified phenomenon as {\it oscillation-induced gain modulation} (OIGM) that we will demonstrate in several examples.

This article is organized as follows: In Section \ref{PPexample}, we will study a first example of a 2D predatory-prey model \cite{holling1965}. In Section \ref{meanOrbitsThm}, we will derive the mean value deviation theorem for general two-dimensional models. In Section \ref{Examples}, we showcase Hopf-induced mean value deviation in several more two-dimensional examples, including the Brusselator model of chemical reactions \cite{belousov1958, prigogine1968}, and a Wilson-Cowan model of neural population activity \cite{wilson1972}.  These examples will demonstrate instances of OIGM in which the cycle mean can dramatically depart from the equilibrium at the bifurcation in functionally significant ways. In Section \ref{nDim}, we will derive the cycle mean deviation theorem for general $n$-dimensional systems ($n>2$) and present a $n=3$-dimensional example.

\section{Example: Predator-Prey Model}
\label{PPexample}

\begin{figure*}
  \centering
  \includegraphics[width=14cm]{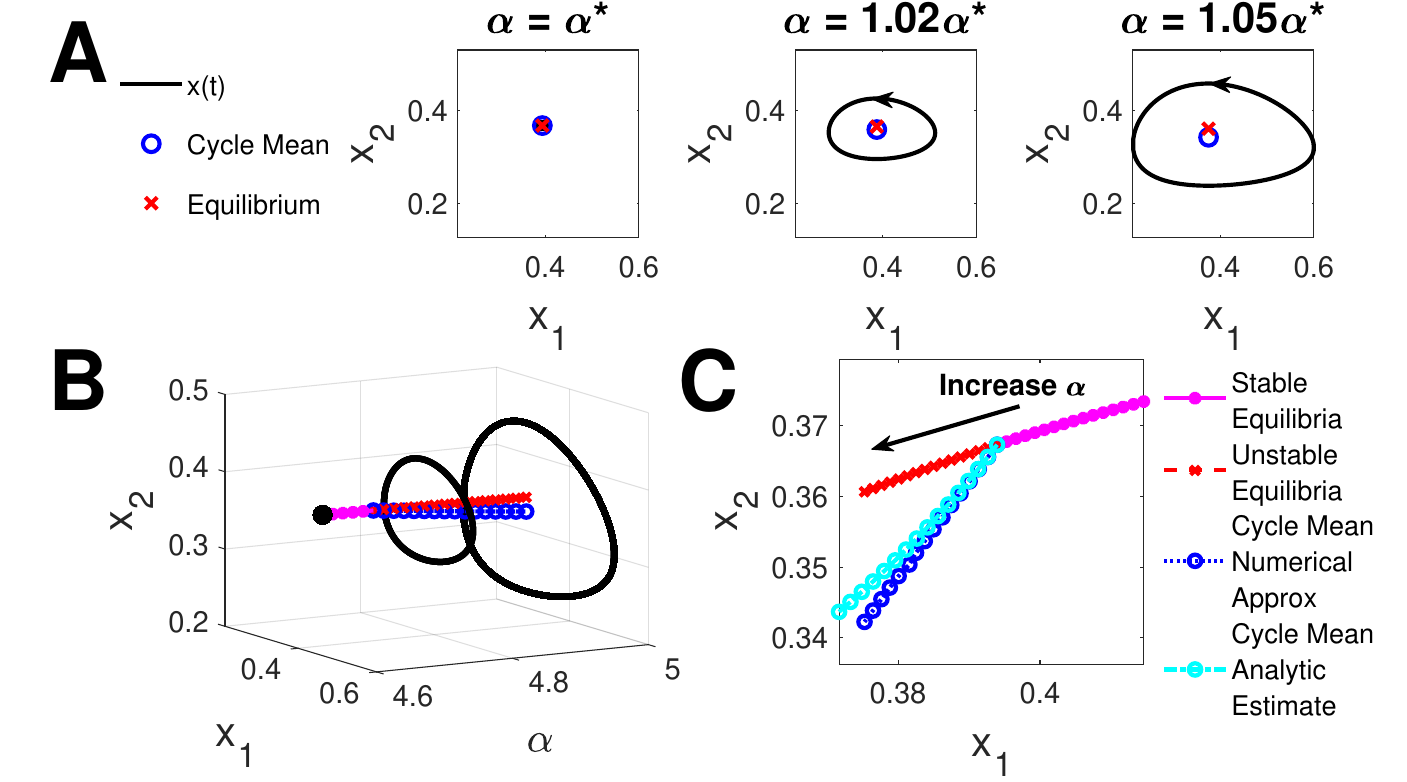}
  \caption{The predator-prey model, consisting of coupled prey ($x_1$) and predator ($x_2$) sizes, illustrates cycle mean deviation from equilibrium. (A) Simulations yield either stable equilibria or stable oscillatory orbits depending on the bifurcation parameter $\alpha$ value relative to the bifurcation point $\alpha^*$. (B) For $\alpha < \alpha^*$ the mean value is the stable equilibrium, but for $\alpha > \alpha^*$ stable cycles emerge (examples shown in black). As $\alpha$ increases beyond $\alpha^*$, the cycle mean (blue) deviates from the unstable equilibria (red). (C) A close-up at the cycle mean and unstable equilibria in phase space, while analytical estimates (cyan) accurately predict numerically computed cycle means (dark blue) near the bifurcation point (see Section \ref{meanOrbitsThm}).}
  \label{Fig1}
\end{figure*}

Consider a predator-prey model \cite{holling1965} in which $x = (x_1, x_2)^T$ is the normalized prey and predator population sizes, respectively, governed by the nonlinear system of differential equations
\begin{align}
\label{predPrey}
\begin{split}
\dot{x}_1 & =  \beta x_1(1-x_1)- c\frac{\alpha x_1}{1+ \alpha x_1} x_2, \\
\dot{x}_2 & = -\delta x_2 + c\frac{\alpha x_1}{1+ \alpha x_1} x_2,
\end{split}
\end{align}
where the first term in the $x_1$-equation (\ref{predPrey}) drives logistic growth of prey, and the second term defines the predation-related decline. The $x_2$-equation's first term defines a constant death-rate for the predator population, while the second term drives predator growth due to predation of prey. The saturating first-order hill-function $c\tfrac{\alpha x_1}{1+ \alpha x_1}$ of the predation term dictates that as prey grow more numerous, there is a protective ``herd" effect which saturates to a maximum $c$ of predation per predator. The $\alpha$-parameter modulates the overall slope of the predation response in which $x_1 = 1/\alpha$ is the half-maximum point $c/2$\footnote{We have chosen parameter symbols consistent with that used in \cite{kuznetsov98}, with the exception of $\alpha \to 1/\alpha$. We inverted $\alpha$ so that it is more easily describable.}. A positive population equilibrium exists provided $c>\delta$: 
 \begin{align}
x_0(\alpha)  =  \Big( \frac{ \delta }{\alpha(c-\delta)}, \quad   \frac{\beta}{\alpha(c-\delta) } \Big[ 1 -  \frac{ \delta }{\alpha(c-\delta)}  \Big] \Big)^T.
\end{align}

To better perform an analysis of the model, it is useful to perform a change of variables from original timescale $t$ to a modified one $\tau$:  $ t = \tau (1+\alpha x_1)$. This change of variables results in a simpler-to-analyze system with polynomial terms:
\begin{align}
\label{predPreyT}
\begin{split}
\dot{x}_1 & = \beta x_1(1-x_1)(1+\alpha x_1)- c\alpha x_1 x_2, \\
\dot{x}_2 & = -\delta x_2(1+\alpha x_1) + c\alpha x_1 x_2.
\end{split}
\end{align}
This modified model (\ref{predPreyT}) is orbitally equivalent to the original, but the temporal evolution of solutions between $x(\tau)$ and $x(t)$ will differ. Therefore, in general, we can expect the temporal averages of oscillatory orbits (\ref{xmean0}) will differ between (\ref{predPreyT}) and (\ref{predPrey}). However, for the pedagogical purposes of this article, we will study the simplified model (\ref{predPreyT}) with parameters $c = 2$, $\delta = 1.3$, and $\beta = 2$, in what follows.

The Jacobian linearization of the system (\ref{predPreyT}) results in eigenvalues with real and imaginary parts
 \begin{align}
 \label{ReevalPredPrey}
 Re(\lambda(\alpha)) &  \equiv \mu(\alpha) = \frac{r\delta(c+\delta)}{2\alpha(c-\delta)}\Big[  \frac{c-\delta}{c+\delta}- \frac{1}{\alpha}\Big], \\
  \label{IMevalPredPrey}
 Im(\lambda(\alpha)) &  \equiv \omega(\alpha) = \frac{ c \sqrt{\beta \delta(c-\delta)}}{(c+\delta)^{\frac{3}{2}}}, 
\end{align}
respectively. A zero real part occurs at $\alpha^* = \frac{c+\delta}{c-\delta}$, and increasing $\alpha$ above $\alpha^*$ induces a loss of equilibrium stability concomitant with the emergence of a small-amplitude oscillation with frequency $\omega>0$ (\ref{IMevalPredPrey}).

Numerical solutions of $x(t)$ at a range of $\alpha$-values illustrate the hallmark feature of Hopf bifurcations of supercritical type: increasing $\alpha$ progressively past $\alpha^*$ elicits increasing amplitude oscillations, starting from zero amplitude, surrounding an unstable equilibrium (Fig. \ref{Fig1}A). Note that the cycle mean (\ref{xmean0}) diverges from the equilibrium location with increasing $\alpha$ (Fig. \ref{Fig1}B). Figure \ref{Fig1}C further illustrates this deviation over a more dense sampling of $\alpha$-values. For $\alpha<\alpha^*$, a stable equilibrium exists and the mean $\langle x \rangle_{\alpha} = x_0(\alpha)$. Then, for $\alpha \geq \alpha^*$, the unstable equilibria continues along on a largely similar trajectory as for sub-threshold $\alpha$ values. However, the numerically computed cycle mean diverges from the unstable equilibria past the bifurcation, producing significantly lower mean predator levels than the equilibrium would predict. Figure \ref{Fig1}C also shows the cycle mean analytic estimate, which for $\alpha$-values near the bifurcation, accurately predicts the numerically computed cycle mean values, but looses accuracy for $\alpha$ further away from the bifurcation. This analytic estimate will be derived below in Section \ref{meanOrbitsThm}. 

\begin{figure*}
  \centering
  \includegraphics[width=14cm]{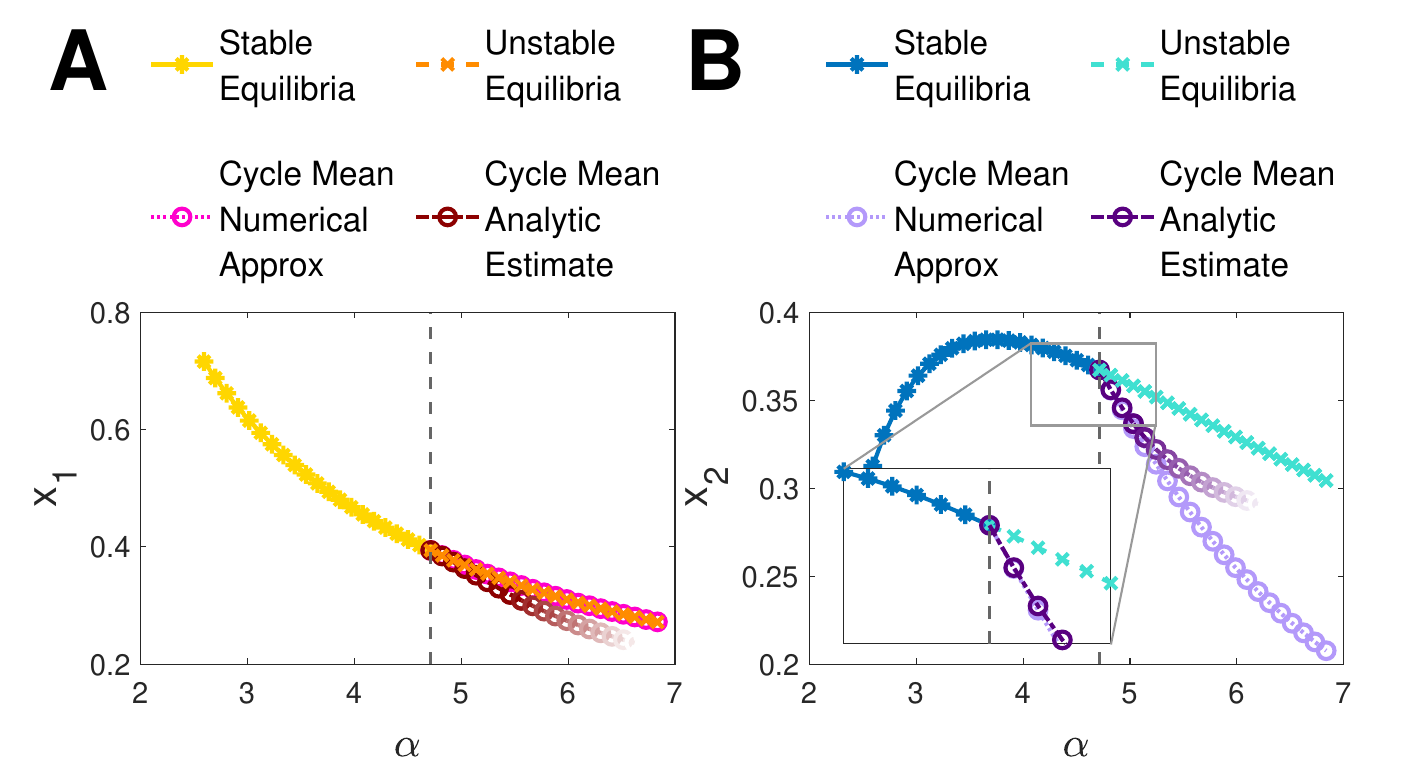}
  \caption{An example of oscillation-induced gain modulation (OIGM). (A) As $\alpha$ increases but remains below $\alpha^*$, the prey ($x_1$) stable equilibria (yellow) decrease. As $\alpha$ increases above $\alpha^*$, the numerically approximated  cycle mean prey values (pink) decrease at the same rate as the unstable prey equilibria (orange x's), while the cycle mean analytic estimate (see Theorem 1) is tangent to the aforementioned numerical estimates near the bifurcation but looses accuracy thereafter. (B) As $\alpha$ increases but remains below $\alpha^*$, predator ($x_2$) equilibria (dark blue) first increase, then decrease. Above $\alpha^*$, the unstable equilibria (cyan) continue on the same trajectory as the stable equilibria; however, the cycle mean numerical approximation, and the analytic estimate (lavender and purple, respectively), show that the mean predator values decrease at a faster rate than the unstable equilibria as a function of $\alpha$. Hence, the oscillation induces an abrupt downward shift in the slope of the mean predators, which we term OIGM.}
  \label{Fig2}
\end{figure*}

Figure \ref{Fig1}C also demonstrates that the cycle mean deviation produced by the Hopf instability dramatically exacerbates the drop in mean $\langle x_2 \rangle $ predator levels relative to drop in prey levels $\langle x_1 \rangle$---i.e., a discontinuity in $ \frac{d\langle x \rangle_{\alpha}}{d\alpha}$ at the bifurcation. This shift in the cycle mean trajectory is our first example of OIGM. Figure \ref{Fig2} shows this OIGM phenomenon as a function of $\alpha$. The mean $x_1$ prey levels shows little change in the gain (i.e., slope) at the bifurcation (Fig. \ref{Fig2}A). However, there is a dramatic shift in the gain of mean predator levels $x_2$ at the bifurcation (Fig. \ref{Fig2}B).

\section{The Hopf Mean Value Theorem for 2-Dimensional Systems}
\label{meanOrbitsThm}
We let $x = (x_1, x_2)^T$, and $\dot{x} = f(x; \alpha) = (f_1, f_2)^T$.  To obtain the $\mathcal{O}(\mu(\alpha)^1)$-estimate of cycle means, we will develop the standard Poincar\'e mapping from $x(t)$ to the normal form equation of a Hopf bifurcation in the complex plane $w(t) \in \mathbb{C}$ \cite{poincare1879,guckenheimer2013,kuznetsov98}. That such a mapping exists is recounted in detail in \cite{kuznetsov98,guckenheimer2013} and originally proven in \cite{arnold2012, arnold2014}. In this article, we have chosen to follow the notation in  \cite{kuznetsov98}. This transformation from 2D real phase space $x$, to the the complex plane enables a description of the limit cycle orbits in polar coordinates with amplitude $r(t)$ and phase angle $\theta(t)$: $w(t) = re^{i\theta}$. We will analyze relevant features of this mapping that enable computing estimates of the mean value (\ref{xmean0}), to lowest order, in terms of $\mu(\alpha)$. 

Without loss of generality, we will re-center the equilibrium to the origin $x \to x-x_0$. Following \cite{kuznetsov98}, we then segment $f$ into linear and non-linear parts $f(x; \alpha) = A(\alpha)x + F(x; \alpha)$, with linear part defined as the Jacobian at $0$:
\begin{align*}
A(\alpha) &  = \begin{bmatrix} a & b \\ c & d \end{bmatrix},
\end{align*}
where $a,b,c,d$ are functions of $\alpha$ given by
\begin{align*}
 a = \frac{\partial f_1}{\partial x_1}(0,\alpha), \quad    b = \frac{\partial f_1}{\partial x_2}(0,\alpha),  \\
 c = \frac{\partial f_2}{\partial x_1}(0,\alpha), \quad  d = \frac{\partial f_2}{\partial x_2}(0,\alpha). 
\end{align*}
The characteristic polynomial of the matrix $A$ is $\lambda^2 - \sigma \lambda + \Delta = 0$, where
\begin{align*}
\sigma(\alpha) & = a+d = Tr(A), \\
\Delta(\alpha) &  = ad-bc = Det(A).
\end{align*}
The Hopf bifurcation condition requires $\sigma(\alpha^*) = 0$, $\sigma'(\alpha^*) \neq 0$, and $\Delta(\alpha^*) = \omega_0^2 > 0 $, which dictates the eigenvalue
$\lambda(\alpha) = \mu(\alpha) + i \omega(\alpha)$ (and its conjugate) is given by
\begin{align*}
 \mu(\alpha)   = \frac{ \sigma(\alpha)}{2} ,  \quad \omega(\alpha)   =  \sqrt{ \Delta - \mu^2}.
\end{align*}

We also define $p(\alpha)$ and $q(\alpha)$ to be the left- and right-hand eigenvectors of $A(\alpha)$: $Aq = \lambda q$, and $p^T A = \bar{\lambda} p^T$, where $\bar{\lambda}$ is the complex conjugate.  Note that $q$, $\bar{q}$, can form an invertible linear map between $\mathbb{C}$ and $\mathbb{R}^2$:
\begin{align}
\label{xFromz}
x = z q(\alpha) + \bar{z}\bar{q}(\alpha). 
\end{align}
To obtain an inverse to (\ref{xFromz}), we must observe the following two facts. Firstly, note that $p$ can always be normalized so that $\langle p, q \rangle = 1$, where $\langle u,v \rangle = \bar{u_1} v_1 + \bar{u_2}v_2$ is the standard inner product. Secondly, note that $p$ and $\bar{q}$ are necessarily orthogonal:
\begin{align*}
 \langle p,  \bar{q}  \rangle &  =  \langle p, \frac{1}{\bar{\lambda}}  A\bar{q}  \rangle  =  \langle A^T p, \frac{ \bar{q} }{\bar{\lambda}}  \rangle =  \langle \bar{\lambda} p, \frac{ \bar{q} }{\bar{\lambda}}  \rangle  = \frac{\lambda}{\bar{\lambda}}  \langle p,  \bar{q}  \rangle.
\end{align*}
From the above, deduce that $ (1-\lambda/\bar{\lambda} ) \langle p,  \bar{q} \rangle  = 0$; and because $ \lambda  \neq \bar{\lambda}$, it implies $\langle p,  \bar{q} \rangle  = 0$.

The inverse of (\ref{xFromz}) then can be obtained by projecting $x$ from (\ref{xFromz}) onto $p$:
\begin{align*}
 \langle p, x \rangle  = \langle p, (z q  + \bar{z}\bar{q})  \rangle   =  \langle p,  z q \rangle +  \langle p,  \bar{z}\bar{q}  \rangle  = z \langle p, q \rangle  + \bar{z}  \langle p,  \bar{q}  \rangle = z. 
\end{align*}

From the projection of $x$ onto $p$, one can obtain the dynamics of $z$:  $\dot{z} = \langle p, \dot{x} \rangle$:
\begin{align}
\label{zDotDefinition}
\begin{split}
\dot{z} & = \langle p,  ( A (z q  + \bar{z}\bar{q}) + F(z q  + \bar{z}\bar{q}) )  \rangle \\
& =  \langle p,  A (z q  + \bar{z}\bar{q}) \rangle  +  \langle p, F(z q  + \bar{z}\bar{q})   \rangle \\
& =  z \langle p,  \lambda q \rangle +  \langle p,  \bar{z}A \bar{q} \rangle +  \langle p, F(z q  + \bar{z}\bar{q})   \rangle \\
& = \lambda z +  \langle p, F(z q  + \bar{z}\bar{q})  \rangle,
\end{split}
\end{align}
which is defined for all $\alpha$, and the above inner product is defined as
\begin{align}
\label{gzzbar}
 g(z,\bar{z}) \equiv  \langle p, F(z q  + \bar{z}\bar{q}) )  \rangle, 
\end{align}
that results in the map of the system into in the complex plane:
\begin{align}
\label{zDotMap}
\dot{z} = \lambda z + g(z,\bar{z}).
\end{align}
To achieve the aforementioned expression of limit cycle orbits in polar coordinates, we must expand $F$ about the equilibrium $x = 0$.
\begin{align}
\label{F3o}
F(x) = \frac{1}{2!} B(x,x) + \frac{1}{3!}C(x,x,x) + \mathcal{O}(||x||^4),
\end{align}
where $F = (F_1, F_2)^T$, and $B$  and $C$ in (\ref{F3o}) are the second- and third-order tensors, respectively. The $B$ tensor is formed from Hessian matrices of each $m$th component $F_m$ of $F$, for $m = 1,2$:
\begin{align*}
B_{m}(x,y)  =  \begin{bmatrix} y_1 & y_2  \end{bmatrix}  \begin{bmatrix} \partial_{11}F_m & \partial_{12}F_m \\  \partial_{21}F_m & \partial_{22}F_m  \end{bmatrix} \begin{bmatrix} x_1 \\ x_2  \end{bmatrix}  = \sum_{j,k = 1}^2( \partial_{jk}F_m)x_k y_j ,
\end{align*}
where $\partial_{jk}$ are 2nd-order partial derivative operators. The $C$ tensor in (\ref{F3o}) is given by
\begin{align*}
C_{m}(x,y,u)  & = u_1 \begin{bmatrix} y_1 & y_2  \end{bmatrix}  \begin{bmatrix} \partial_{111}F_m & \partial_{112}F_m \\  \partial_{121}F_m & \partial_{122}F_m  \end{bmatrix} \begin{bmatrix} x_1 \\ x_2  \end{bmatrix} +  u_2 \begin{bmatrix} y_1 & y_2  \end{bmatrix}  \begin{bmatrix} \partial_{211}F_m & \partial_{212}F_m \\  \partial_{221}F_m & \partial_{222}F_m  \end{bmatrix} \begin{bmatrix} x_1 \\ x_2  \end{bmatrix}  \\
& =  \sum_{j,k,\ell = 1}^2 (\partial_{\ell jk}F_m )x_k y_j  u_{\ell},
\end{align*}
where $\partial_{jk\ell}$ are third-order partial derivative operators. 

Evaluating the low-order expansion (\ref{F3o}) that includes $B(x,x)$ and $C(x,x,x)$ in terms of $z$ using the linear map (\ref{xFromz}) results in 
\begin{align}
\label{Bz}
\begin{split}
B( z q  + \bar{z}\bar{q} , z q  + \bar{z}\bar{q})  = z^2 B(q,q) + 2z\bar{z}B(q,\bar{q})+ \bar{z}^2 B(\bar{q},\bar{q}),
\end{split}
\end{align}
while for $C$, this results in
\begin{align}
\label{Cz}
\begin{split}
C( z q  + \bar{z}\bar{q} , z q  + \bar{z}\bar{q}, z q  + \bar{z}\bar{q}) & = 
C(q,q,q) z^3  + 3C(q,q,\bar{q})z^2\bar{z}+3C(q,\bar{q},\bar{q}) z\bar{z}^2 +  C(\bar{q},\bar{q},\bar{q})\bar{z}^3 .
\end{split}
\end{align}

The above expression of $F$ in terms of the multilinear functions of $z$, $\bar{z}$ (\ref{Bz}-\ref{Cz}) enable the expansion of (\ref{gzzbar}) in terms of powers of $z$, $\bar{z}$:
\begin{align}
\label{projExpansion}
\begin{split}
\langle p, F(z q  + \bar{z}\bar{q}) ) \rangle &  =    \frac{1}{2}( \langle p, B(q,q) \rangle z^2   + \langle p, B(q,\bar{q}) \rangle 2z\bar{z}   +   \langle p, B(\bar{q},\bar{q})  \rangle\bar{z}^2  ) \ldots  \\
&  + \frac{1}{6}( z^3  \langle p, C(q,q,q)  \rangle + 3  \langle p, C(q,q,\bar{q})  \rangle z^2\bar{z}  \ldots \\
&  +3\langle p, C(q,\bar{q},\bar{q})  \rangle z\bar{z}^2  +  \langle p, C(\bar{q},\bar{q},\bar{q})  \rangle  \bar{z}^3  )  + \mathcal{O}(|z|^4),
\end{split}
\end{align}
where the inner products in (\ref{projExpansion}) define a set of coefficients $g_{jk}$: 
\begin{align}
\label{gcoef}
\begin{split}
g_{20} & =  \langle p, B(q,q) \rangle, \\
g_{11} & = \langle p, B(q,\bar{q})  \rangle,  \\
g_{02} & = \langle p, B(\bar{q},\bar{q})  \rangle, \\
g_{30} & =  \langle p, C(q,q,q)  \rangle, \\
g_{21} & =  \langle p, C(q,q,\bar{q})  \rangle, \\
g_{12} & =  \langle p, C(q,\bar{q},\bar{q})  \rangle, \\
g_{03} & =  \langle p, C(\bar{q},\bar{q},\bar{q})  \rangle.
\end{split}
\end{align}
Note that $g_{11}$ is real-valued because $\overline{ B(q,\bar{q}) }  = B(\bar{q},q)$. The above coefficients (\ref{gcoef})  define the DE for $z$ up to third-order:
\begin{align}
\label{zdotPoly}
\begin{split}
\dot{z} &  = \lambda z + g(z,\bar{z}) \\
     & =  \lambda z + \sum_{k+j \geq 2}^{3} \frac{1}{k! j!}g_{jk}(\alpha) z^j \bar{z}^k + \mathcal{O}(|z|^4) 
\end{split} 
\end{align}
The DE (\ref{zdotPoly}) defines, to lowest order, the solutions $z(t)$ near the origin. By itself, the description in the $z$-variable is no more useful than the original $x$. To make progress, we use an additional near-identity transformation to a new complex variable $w$ that is governed by ay simplified Poincar\'e normal form DE, whose orbits are homeomorphic to those in $z$. Consider the near-identity mapping
  \begin{align}
  \label{zFromw}
z  &  =  w + \frac{h_{20}}{2}w^2 + h_{11}w\bar{w} + \frac{h_{02}}{2}\bar{w}^2   
+ \frac{h_{30}}{6}w^3 + h_{21}w^2\bar{w} + \frac{h_{12}}{2}w\bar{w}^2 +  \frac{h_{30}}{6}w^3,
\end{align}
and its inversion, to lowest order
\begin{align}
\label{wfromz}
\begin{split}
  w  & = z- (\frac{h_{02}}{2}w^2 + h_{11}w\bar{w} + \frac{h_{20}}{2}\bar{w}^2 ) + \mathcal{O}(|w|^3)  \\
  & = z - \frac{h_{02}}{2}z^2 - h_{11}z\bar{z} - \frac{h_{20}}{2}\bar{z}^2 + \mathcal{O}(|z|^3). 
  \end{split}
\end{align}
Following \cite{arnold2012,arnold2014,  kuznetsov98}, we differentiate the above equation (\ref{wfromz}) to obtain $\dot{w}$ in terms of $z$, $\bar{z}$, $\dot{z}$, and $\dot{\bar{z}}$. After substituting expressions involving the $g_{jk}$, and combining like terms involving $h_{jk}$, the $h_{jk}$ can be chosen so that $\dot{w}$ depends, to lowest order, solely on $\lambda w$, and a single third order term
 \begin{align}
  \label{wdot}
 \dot{w} = \lambda w +c_1 w^2 \bar{w} + \mathcal{O}(|w|^4).
  \end{align}
That is, all other second- and third-order terms vanish (see \cite{kuznetsov98}) except for a single third order term in $w^2\bar{w}$, with coefficient $c_1$ \cite{Bautin1949} given by
     \begin{align}
  \label{c1}
c_1 = \frac{2\lambda + \bar{\lambda}}{2|\lambda |^2} g_{20}g_{11} + \frac{1}{\lambda} |g_{11}|^2 + \frac{1}{2(2\lambda - \bar{\lambda})} |g_{02}|^2 + \frac{1}{2}g_{21}.
  \end{align}
  Note that we assume that $Re(c_1) \neq 0$, which defines a codimension one, simple Hopf point with a single limit cycle. The quadratic values $h_{jk}$ that achieve this cancellation in (\ref{wdot}) are found to be
\begin{align}
\label{hquadg}
  h_{20} =   \frac{g_{20}}{\lambda}, \quad  h_{11} = \frac{g_{11}}{\bar{\lambda} }, \quad  h_{02} = \frac{g_{02} }{2\bar{\lambda}- \lambda }. 
  \end{align}
The third order terms can be cancelled similarly, with the exception of the $w^2 \bar{w}$ term \cite{kuznetsov98}:
\begin{align}
\label{hcubic}
  h_{30} =  \frac{g_{30}}{2\lambda}, \quad  h_{21} = 0, \quad  h_{12} =  \frac{g_{12}}{2\bar{\lambda} }, \quad  h_{03} = \frac{g_{03} }{3\bar{\lambda}- \lambda }. 
\end{align}
Note that $h_{21} = 0$ is chosen because to cancel the $z^2 \bar{z}$-term involves zeroing the equation $g_{21} - 2\mu(\alpha) h_{21} $ which would cause $h_{21}$ to otherwise diverge as $\alpha \to \alpha^*$ (see \cite{kuznetsov98}). 

 As we shall show, the most relevant of the $h_{jk}$ for the purposes of finding a mean deviation from equilibrium $\langle x \rangle_{\alpha} -x_0(\alpha)$ is the $h_{11}$-term: $h_{11} = g_{11}/\bar{\lambda}$.

 The order-four and higher terms can be neglected near the origin, where we seek small amplitude circular orbit solutions of the form $w(t) = r_w e^{i \omega_w t}$ in (\ref{wdot}), where we label $r_w$ and $\omega_w$ as the radius and frequency, respectively, in the complex $w$-plane denoted by ``w". Substituting $w(t) = r_we^{i \omega_w t}$ into (\ref{wdot}), we find:
  \begin{align}
\label{r_wDef}
 r_{w} & = \sqrt{ \frac{ -\mu}{Re(c_1)}}, \\
 \label{omega_wDef}
 \omega_{w} & = \omega - \frac{Im(c_{1})}{Re(c_{1})}\mu.  
 \end{align}
 Substituting the limit cycle solution $w(t)$ (\ref{r_wDef}-\ref{omega_wDef}) into (\ref{zFromw}), one obtains a limit cycle solution in the $z(t)$ plane.  With the maps (\ref{zFromw}) and (\ref{xFromz}), we can then project the circular $w(t)= r_w e^{i \omega_w t}$ solution into the $x$-plane to get an approximation of the $x(t)$ limit cycle. This mapping is a suitable approximation when $\alpha$ is near the bifurcation point. The cycle mean of $x$ can be computed term-by-term through the mapping (\ref{zFromw}). As we will show below, most terms in the mapping integrate to zero, except one term, greatly simplifying the calculation. The result is thus:

{\bf theorem}: The Hopf Mean Deviation Theorem for 2-Dimensional Systems  

Consider a smooth 2D system $\dot{x} = f(x; \alpha)$ undergoing a codimension one Hopf bifurcation at $\alpha^*$, inducing an oscillatory solution $x(t)$ for $\alpha \in \tilde{I}_{\alpha^*}$ surrounding an equilibrium $x_0(\alpha)$. The oscillatory mean $\langle x \rangle_{\alpha} = \frac{1}{T} \int_0^T x(t) dt$ deviates from the equilibrium according to
  \begin{align}
 \label{mainResultT}
 \langle x \rangle_{\alpha} -  x_0(\alpha)  =  \begin{cases}  K(\alpha) \mu(\alpha)  + \mathcal{O}(\mu(\alpha)^2), &  \alpha \in \tilde{I}_{\alpha^*} \\ 0,  &  \alpha \in I_{\alpha^*}  \end{cases}.
\end{align}
The vector quantity $K$ in (\ref{mainResultT}) is given by
\begin{align}
   \label{Kfinal}
  K =   - Re \Big( 2 \frac{g_{11} }{\bar{\lambda} Re(c_{1})}\Big) Re(q)  + Im \Big( 2 \frac{g_{11} }{ \bar{\lambda}  Re(c_{1}) }\Big) Im( q)  ,
  \end{align}
  where $g_{11}$, and $c_{1}$ are defined above in (\ref{gcoef}) and (\ref{c1}).

{\bf proof:}

The proof involves temporally averaging the $z(t)$ expression defined in terms of $w(t)$ (\ref{zFromw}). Note, however, when $w(t) = r_w e^{i \omega_w t}$ the integrals compute to zero for all cases except $j = k$:
 \begin{align}
 \frac{1}{T} \int_0^T w^j(t)\bar{w}^k(t) dt =   \begin{cases} 0 , &  j \neq k \\ r_w^{j+k}, & j = k \end{cases}. 
 \end{align}
That is, for lowest orders below 3, only the quadratic term consisting of factors $w^1$ and $\bar{w}^1$, results in a non-zero result. That is, $\frac{1}{T} \int_0^T w(t) \bar{w}(t) dt= r_w^2. $
  Therefore, averaging (\ref{zFromw}) over the period $T$ results in 
   \begin{align}
   \label{zmean}
\frac{1}{T}\int_0^T z(t) dt = h_{11}r_w^2  =    -\frac{g_{11}}{\bar{\lambda} Re(c_1)} \mu +  \mathcal{O}(\mu^2).
  \end{align}
  Using (\ref{zmean}), we average over $x = z q + \bar{z}\bar{q}$ (\ref{xFromz}) to obtain the desired mean of $x(t)$ from its equilibrium:
   \begin{align}
   \label{kProof}
\frac{1}{T}\int_0^T x(t) dt - x_0(\alpha) =  K(\alpha) \mu(\alpha) + \mathcal{O}(\mu^2),
  \end{align}
 where the vector quantity $K(\alpha)$ is given by
    \begin{align}
   \label{Kfinal}
K(\alpha) =   - Re \Big( 2 \frac{g_{11}}{ \bar{\lambda} Re(c_{1}) } \Big) Re(q)  + Im \Big( 2 \frac{g_{11}}{ \bar{\lambda} Re(c_{1}) } \Big) Im(q)  .
  \end{align}

$\square$

Remark: The real-valued coefficient $g_{11}$ of (\ref{gcoef}) is the key factor in (\ref{Kfinal}) whose value dictates the  size of the mean deviation (\ref{mainResultT}). However, $g_{11}$ is not sole determiner of other reported features of Hopf bifurcations. In particular, the denominator of (\ref{Kfinal}) contains tensor-derived value $ Re(c_1(\alpha))$, termed the first Lyapunov coefficient $\ell_1$, which is a familiar term in the study of Hopf bifurcations, and is nonzero in the case of a codimension-one Hopf bifurcation. The factor $c_1$ (\ref{c1}) contains many tensor derived coefficients including $g_{21}$, $g_{20}$, as well as $g_{11}$. Hence, $g_{11}$ is the most prominent coefficient determining mean deviation from equilibrium.

\section{Additional Examples}
\label{Examples}

\subsection{Analysis of predator-prey model}
\label{predPreyAnalysis}
Returning to our original predator-prey model example, it is computed that 
\begin{align}
\label{predPreyGs}
\begin{split}
g_{20} &= \frac{ c\delta(c^2-\delta^2-\beta \delta ) + i\omega c(c+\delta)^2 }{c+\delta}, \\
g_{11} & = -\frac{\beta c \delta^2}{c+\delta}, \\
g_{21} &  = -3rc^2 \delta^2.
\end{split}
\end{align}
Plugging in the eigenvalues (\ref{ReevalPredPrey},\ref{IMevalPredPrey}) and coefficients (\ref{predPreyGs}) into (\ref{c1}) leads to the evaluation of  (\ref{Kfinal}). The expression for $K(\alpha)\mu(\alpha)$ is large and we won't write it fully here, but the calculated analytical estimate well-approximates those computed from numerical simulations near the bifurcation (see Figs. \ref{Fig1}-\ref{Fig2}).

\subsection{Analysis of the Brusselator model}
\label{BrusModel}

The Brusselator equations model chemical reactions that can exhibit limit cycle solutions \cite{prigogine1968}. Let $x_1$ be the quantity of a chemical reagent, and $x_2$ is the quantity of some chemical product made from the reagent, governed by the differential equations (DEs)
\begin{align}
\begin{split}
\dot{x}_1 & =  A - (\alpha+1)x_1 +  x_1^2 x_2,  \\
\dot{x}_2 & = \alpha x_1 - x_1^2 x_2.
\end{split}
\end{align}
The parameter $A$ is an input growth rate of the $x_1$-species, while $\alpha$ is the sets the linear rate of conversion of $x_1$ to $x_2$. The nonlinear last term in each of the above equations produces negative feedback if the product or reagent grows too large. The equilibrium is found to be $(x_1^*, x_2^*) = (A, A/\alpha)$. We take $\alpha$ to be the bifurcation parameter. The equilibrium becomes unstable when $\alpha \geq 1 + A^2$. The bifurcation point $\alpha^* = 1+A^2$ then depends on $A$. 

\begin{figure*}
  \centering
  \includegraphics[width=14cm]{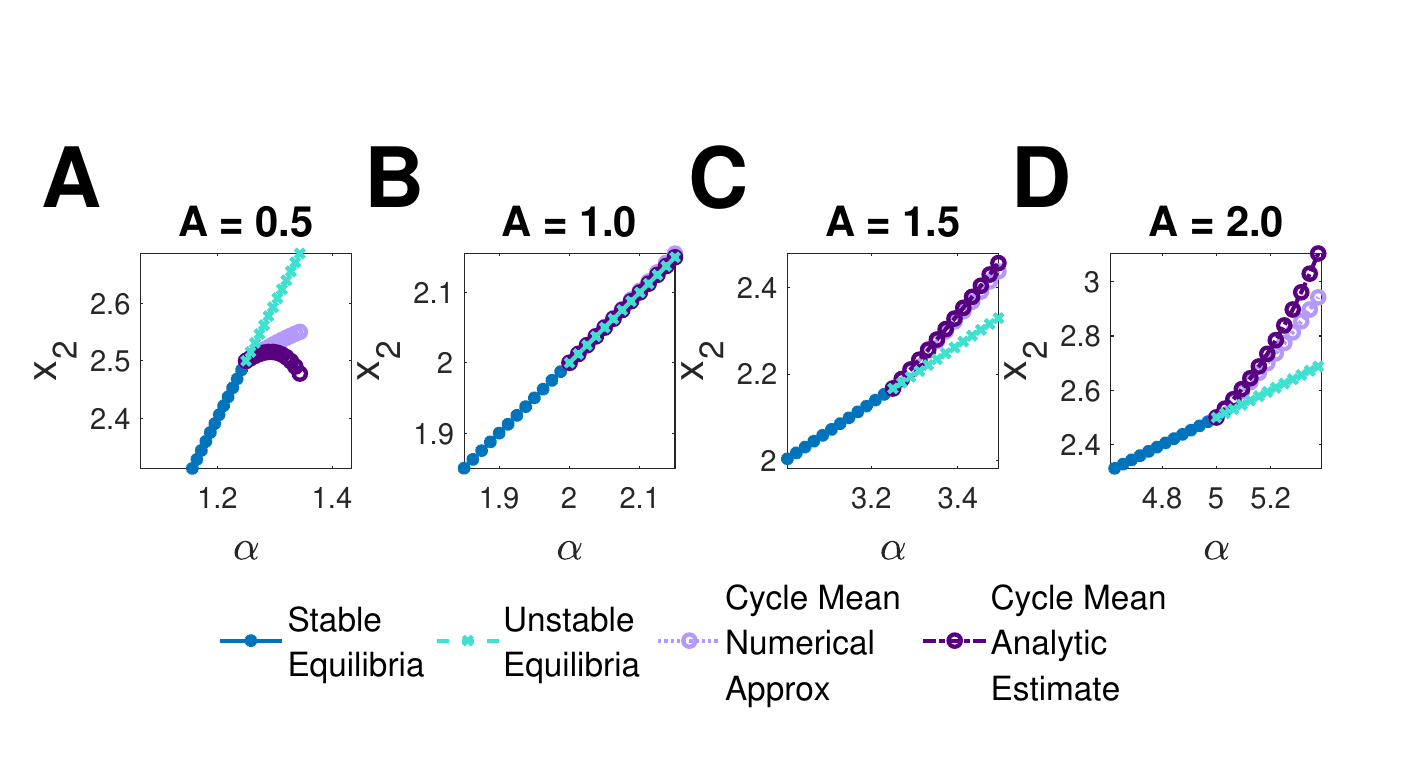}
  \caption{The Brusselator model with the mean of $x_2$ as a function of $\alpha$, over various $A$-values. (A) For $A < 1$, both the analytically estimated and numerically approximated cycle means (lavender and purple circles) modulate below the unstable equilibrium (cyan). (B) For $A = 1$, the cycle mean coincides with the unstable equilibrium because $g_{11} = 0$. (C,D) For $A > 1$, both the analytically estimated and numerically approximated $x_2$ mean deviates are above its unstable equilibrium value. }
  \label{Fig3}
\end{figure*}

Figure \ref{Fig3} shows the stable and unstable equilibria in phase space, as well as the numerical and analytic cycle mean estimates for several distinct $A$-values. The relevant Poincar\'e mapping coefficients are
\begin{align}
\label{gBrussel}
\begin{split}
g_{20} = A-i, \quad g_{11} = \frac{(A-i)(A^2-1)}{A^2 + 1}, \quad g_{21} = \frac{A(3A-i)}{A^2+1}.
\end{split}
\end{align}
Values of $A$ below unity induce a reduction in gain (i.e. slope) when Hopf bifurcation is induced (Fig. \ref{Fig3}A). Additionally, note that for $A=1$, equation (\ref{gBrussel}) dictates that $g_{11}=0$, so that $K = 0$, meaning there is zero deviation of the cycle mean from the equilibrium (Fig. \ref{Fig3}B). Moreover, when $A$ is increased to values progressively greater than unity, there is an increase in gain when a Hopf bifurcation is induced (Fig. \ref{Fig3}C, D). 

This example underscores the determinative nature of the $g_{11}$ coefficient on the mean deviation.  This example also demonstrates that the gain of mean cycle deviation with respect to the bifurcation parameter ($\alpha$) can be dramatically manipulated by other model parameters ($A$) and can serve to either increase or decrease the mean response to $\alpha$ due to Hopf-induced oscillations.

\subsection{Analysis of the Wilson-Cowan model}
\label{WCModel}
The Wilson-Cowan model describes the synaptic activity of populations of excitatory (e) and inhibitory (i) neurons by variables $u_1$ and $u_2$, respectively \cite{wilson1972}: 
\begin{align}
\label{WCmodel}
\begin{split}
\tau_1\dot{u}_1 &= -u_1 + f(w_{ee} u_1 - w_{ei} u_2 + I), \\
\tau_2\dot{u}_2 &= -u_2 + f(w_{ie} u_1 - w_{ii} u_2 + I). \\
\end{split} 
\end{align}
The $\tau_j$ are timescales that set the rate of change of the synaptic variables to track the current firing rate $f$ of the respective e- and i-populations (\ref{WCmodel}). The firing rate function $f$ is defined by
\begin{align}
\label{sigmoid}
\begin{split}
f(x) &= (1 + e^{ \frac{-x}{\beta}})^{-1},
\end{split} 
\end{align}
and possesses a sigmoid shape in which $\beta$ varies the steepness of the sigmoid function.  The firing rate of each population is the function of the weighted synaptic inputs, in which $w_{jk}$ are the weights within and between e- and i- cells, and an external input current $I$. The bifurcation parameter in this example is the input current $I$ and supercritical Hopf instabilities have been observed in previous studies \cite{powanwe21}. In this example, we have chosen values $\tau_1 = 4$, $\tau_2 = 12$, $\beta = 0.1$,  and $w_{ee} = 3.6$, $w_{ei} = 8$, $w_{ie} = 4$, $w_{ii} = 8.8$. 

\begin{figure*}
  \centering
  \includegraphics[width=14cm]{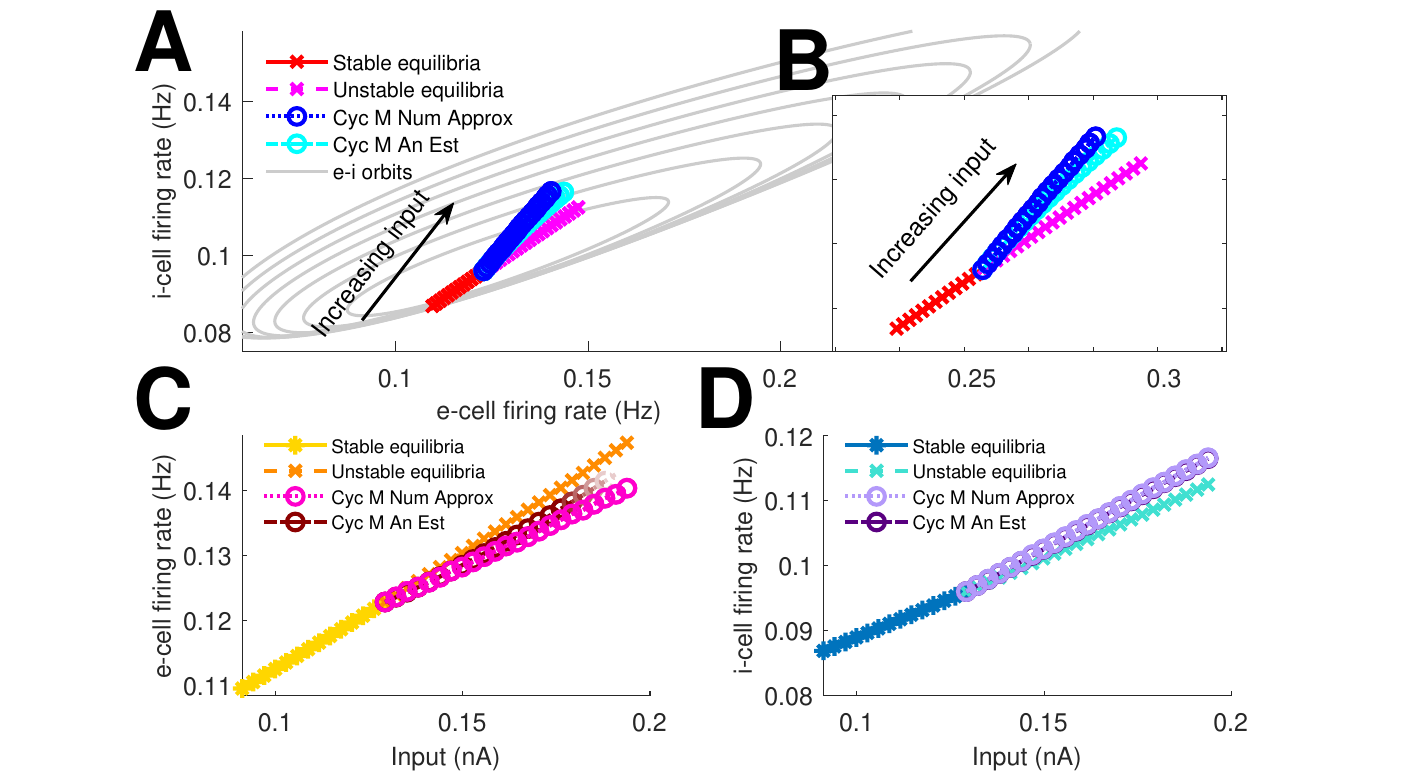}
  \caption{The Wilson-Cowan model exhibits OGIM in both the e- and i-cell populations, but in opposite directions. (A) As input increases beyond the bifurcation point, oscillations emerge with progressively larger amplitudes, show in phase space. The cycle mean also exhibits a trajectory shift at the bifurcation. (B) a close-up of panel A, showing oscillation-induced cycle mean trajectory shifts for both numerically and analytically computed cycle mean estimates. (C) e-cell exhibits oscillation-induced reduction in cycle mean gain as a function of $I$. (D) The opposite occurs for i-cells.}
  \label{Fig4}
\end{figure*}

With the above parameters, sufficiently low values of input $I$ elicit a stable equilibrium in which e- and i-cell activity find a balance point. As $I$ increases, the location of the stable equilibrium location shifts upward for both e- and i-cells, as shown in Figure \ref{Fig4}A (see also inset B). Increasing $I$ further past a threshold $I^*$ destabilizes the equilibrium, and a stable limit cycle emerges with progressively larger oscillation amplitudes  (Fig. \ref{Fig4}A). The onset of oscillations coincides with a shift in the mean activity trajectory. Past the bifurcation, the e-cell mean activity exhibits reduction in gain (Fig. \ref{Fig4}C), while i-cells show an uptick in gain (Fig. \ref{Fig4}D) as a function of $I$. This oscillation-induced role-reversal between e- and i-cell mean activity may have consequences to any postsynaptic target that is sensitive to the average activity of  e- and/or i-cell outputs.

\section{The Hopf Mean Value Theorem for $n$-Dimensional Systems}
\label{nDim}

The $\mathcal{O}(\mu)$ cycle mean estimate can also be computed for $n$-dimensional systems. Following Section \ref{meanOrbitsThm}, we begin with the general equation
\begin{align}
\dot{x} = f(x; \alpha) = Ax + F(x),
\label{xdotndim}
\end{align}
where $x \in \mathbb{R}^n$. We assume there exists an equilibrium that we translate to the origin $x_0(\alpha) = 0$ that exists over a range of $\alpha$ values. We define the linearized part of $f$, expanded about $x_0(\alpha)$, to be the $n\times n$ Jacobian matrix $Df = A$. Assume the linearization has a single pair of complex conjugate eigenvalues $\lambda(\alpha) = \mu(\alpha) \pm i \omega(\alpha)$, such that $\mu(\alpha^*) = 0$. We also assume $\omega(\alpha) >0$, and $\mu'(\alpha^*) \neq 0$. The pair of eigenvalues are associated with left eigenvectors $p$, $\bar{p}$ and right eigenvectors $q$, $\bar{q}$, and a normalization condition, respectively:
\begin{align}
Aq = \lambda q, \quad  p^TA  = \bar{\lambda} p^T, \quad \langle p, q \rangle = 1.
\end{align}
We will also assume all other eigenvalues have negative real part for all $\alpha$ in a domain near the bifurcation. Just as we have demonstrated in Section \ref{meanOrbitsThm}, note that $\langle p, q \rangle = 1$, and $\langle p, \bar{q} \rangle = 0$.

When $\alpha = \alpha^*$, the center eigenspace $T^c$ corresponds with eigenvalues $\lambda = \pm i \omega_0$, and moreover, as $\alpha$  passes the bifurcation point and enters $\tilde{I}_{\alpha^*}$, the center eigenspace becomes the unstable eigenspace $T^u$. These eigenspaces are both two-dimensional and spanned by $ \{Re(q), Im(q)\}$. The stable $(n-2)$-dimensional eigenspace $T^{s}$, is associated with all other eigenvalues. We can decompose $x \in \mathbb{R}^n$ as 
\begin{align}
x = zq + \bar{z} \bar{q} + y
\label{xEquationNDim}
\end{align}
where $y \in \mathbb{R}^n$ is the component of $x$ in $T^{s}$, with the condition that $\langle q, y \rangle = 0$. Also, we define $z \in \mathbb{C}$ being a complex coordinate that locates $x$ in the $T^{c}$ or $T^{u}$ space. We express $z$ and $y$ in terms of $x$ by projecting onto $p$:
\begin{align}
\begin{split}
z  & = \langle p,x \rangle , \\
y & = x - \langle p,x \rangle q - \langle \bar{p}, x \rangle \bar{q}.
\end{split}
\end{align}

Using the original equation \eqref{xdotndim} we can derive the following system of DEs:
\begin{align}
\label{zdotn}
\begin{split}
\dot{z} &= \lambda z + \langle p,F(zq + \bar{z}\bar{q} + y) \rangle, 
\end{split} \\
\label{ydotn}
\begin{split} 
\dot{y} &=  Ay +F(zq + \bar{z} \bar{q} + y) - \langle p, F(z q + \bar{z} \bar{q} + y) \rangle q  - \langle \bar{p}, F(zq + \bar{z} \bar{q} + y) \rangle \bar{q}.
\end{split}
\end{align}

Here $\dot{z}$ (\ref{zdotn}) is similar to its two-dimensional counterpart \eqref{zDotDefinition}, but with the inclusion $y$.
In equations (\ref{zdotn}) and (\ref{ydotn}), we expand $F$ into lowest order multilinear components $F(x) = \frac{1}{2}B(x,x) + \frac{1}{6}C(x,x,x)$ in powers of $z$, $\bar{z}$, as well as including terms with $y$, to yield
\begin{align}
\label{zdotnG}
\begin{split}
\dot{z} &= \lambda  z + \tfrac{1}{2} g_{20} z^2 + g_{11} z \bar{z} + \tfrac{1}{2} g_{02} \bar{z}^2 + \tfrac{1}{2} g_{21} z^2 \bar{z} \ldots \\
 & \quad +  \langle p, B(q,y) \rangle z + \langle p, B(\bar{q},y) \rangle \bar{z} + \cdots, 
 \end{split} 
 \end{align}
 and
 \begin{align}
 \label{ydotnH}
 \begin{split}
\dot{y} &= Ay + \tfrac{1}{2} H_{20} z^2 + H_{11} z \bar{z} + \tfrac{1}{2} H_{02} \bar{z}^2 + \cdots.
\end{split}
\end{align}

The second- $B(x,y)$ and third-order tensors $C(x,y,z)$ of $F$ define the relevant $g_{jk}$ coefficients in $\mathbb{C}$:
\begin{align}
g_{20} = \langle p, B(q,q) \rangle, \quad g_{11} = \langle p, B(q,\bar{q}) \rangle, \quad g_{02} = \langle p, B(\bar{q},\bar{q}) \rangle, \quad g_{21} = \langle p, C(q,q,\bar{q}) \rangle.
\end{align}
The $H_{jk} \in \mathbb{C}^n$ in (\ref{ydotnH}) are defined in terms of $B$ within (\ref{ydotn}) as
\begin{align}
\label{HjkDefn}
\begin{split}
H_{20} =  B(q,q) -  \langle p, B(q,q) \rangle q -  \langle \bar{p}, B(q,q) \rangle \bar{q}, \\
H_{11} =  B(q,\bar{q}) - \langle p, B(q,\bar{q}) \rangle q -  \langle \bar{p}, B(q,\bar{q}) \rangle \bar{q}, \\
H_{02} = B(\bar{q},\bar{q}) - \langle p, B(\bar{q},\bar{q}) \rangle q -  \langle \bar{p}, B(\bar{q},\bar{q}) \rangle \bar{q}.
\end{split}
\end{align}
Note that that $\overline{H_{02}} = H_{20}$ in  (\ref{HjkDefn}), and moreover $H_{11} \in \mathbb{R}^n$, consistent with the fact that $y \in  \mathbb{R}^n$. 

Note also that $x = zq + \bar{z}\bar{q}$ lies in $T^{c}$ and $T^u$ when $\alpha$ is at or past the bifurcation, respectively. However, we must obtain the dynamics of $x$ on the local center $W^c_{loc}(x_0)$ and unstable $W^u_{loc}(x_0)$ invariant manifolds that are tangent to the $T^c$ and $T^u$ subspaces, respectively. To move $x$ into the these invariant manifolds, we must add the right amount of correction $y \in T^s$:
\begin{align}
\label{Wlock}
x = zq + \bar{z}\bar{q} + y.
\end{align}
The correction is written to lowest order in terms of the $z$ and $\bar{z}$ coordinates:
\begin{align}
\label{yWc}
y = V(z, \bar{z}) = \tfrac{1}{2} \eta_{20} z^2 + \eta_{11} z \bar{z} + \tfrac{1}{2} \eta_{02} \bar{z}^2 + O(|z|^3),
\end{align}
where $y$ is orthogonal to $T^{c}$, $T^{u}$, so that $\langle p, y \rangle=0$. Therefore, the vectors $\eta_{ij} \in \mathbb{C}^n$ must also satisfy the orthogonality constraint $\langle p, \eta_{jk} \rangle = 0$. 
To find $\eta_{jk}$, we first equate the known expression of $\dot{y}$ in (\ref{ydotnH}) to the derivative of (\ref{yWc}) 
\begin{align}
\label{dV}
\dot{y} =   \eta_{20} z\dot{z} + \eta_{11} (\dot{z} \bar{z}  + z \dot{\bar{z}})+ \eta_{02} \bar{z}\dot{\bar{z}} + \ldots  \\
\label{dy}
= Ay + \tfrac{1}{2} H_{20} z^2 + H_{11} z \bar{z} +  \tfrac{1}{2} H_{02} \bar{z}^2 + \ldots. 
\end{align}
The above equivalence between (\ref{dV}) and (\ref{dy}) relates the $\eta_{jk}$ to the $H_{jk}$ via the following linear equations and solutions:
\begin{align}
\label{Heqns}
\begin{split}
(2 \lambda I_n - A) \eta_{20}  = H_{20} & \implies \eta_{20} = (2 \lambda I_n - A)^{-1}H_{20},   \\
(2 Re(\lambda) I_n - A )\eta_{11} = H_{11} & \implies  \eta_{11} = (2 Re(\lambda) I_n - A )^{-1}H_{11} , \\
(2 \bar{\lambda} I_n - A) \eta_{02} = H_{02} & \implies  \eta_{02} = (2 \bar{\lambda} I_n - A)^{-1}H_{02},
\end{split}
\end{align}
where $I_n$ is the $n \times n$ identity matrix, and the left-hand-side matrices in (\ref{Heqns}) are invertible, owing to $2\lambda$ and $2Re(\lambda)$ not being eigenvalues of $A$. Thus, the equations (\ref{Heqns}) have unique solutions. The $\eta_{jk}$ solutions in (\ref{Heqns}, right hand side) define, to lowest order, the $y$-correction (\ref{yWc}) that can be substituted into (\ref{zdotnG}), obtaining the lowest-order approximation of the dynamics of $x(t)$ on the center- or unstable manifolds in terms of the coordinate $z$: 
 \begin{align}
\label{zdotnGnoY}
\begin{split}
\dot{z} &= \lambda z + \tfrac{1}{2} g_{20} z^2 + g_{11} z \bar{z} + \tfrac{1}{2} g_{02} \bar{z}^2 + \tfrac{1}{2} g_{21} z^2 \bar{z} \ldots \\
 & \quad \quad +  \langle p, B(q, V(z,\bar{z} ) \rangle z + \langle p, B(\bar{q},V(z,\bar{z} )) \rangle \bar{z} + \cdots, 
 \end{split} 
\end{align}
Substituting $V(z,\bar{z})$ (\ref{yWc}) into (\ref{zdotnGnoY}), allows us to collect up like terms to lowest order, to obtain
 \begin{align}
\label{zdotnGw}
 \begin{split}
\dot{z} &\quad =   \lambda z + \tfrac{1}{2} g_{20} z^2 + g_{11} z \bar{z} + \tfrac{1}{2} g_{02} \bar{z}^2  + \ldots \\
 & \quad \quad + (  \tfrac{1}{2} g_{21} +  \langle p, B(q,  \eta_{11} ) \rangle  + \tfrac{1}{2}   \langle p, B(q,  \eta_{20} ) \rangle  )z^2 \bar{z} + \ldots.
 \end{split}
\end{align}
In the above (\ref{zdotnGw}) the third order term is modified due to the inclusion of the $V(z,\bar{z})$-related influence, which, for simplicity, we denote as $\tilde{g}_{21}$:
 \begin{align}
\label{zdotnG21}
 \begin{split}
 \tilde{g}_{21} =   g_{21} + 2 \langle p, B(q,  \eta_{11} ) \rangle  + \langle p, B(q,  \eta_{20} ) \rangle. 
  \end{split}
\end{align}

The above dynamics of $z(t)$ (\ref{zdotnGw}) approximates that of $x(t)$ on the invariant 2D center- and unstable manifolds embedded in $\mathbb{R}^n$ near the equilibrium. As in the $n = 2$ case, we express $z$ through a near identity mapping in terms of a new variable $w \in \mathbb{C}$, with the purposes simplifying the dynamics of $w$ to only essential terms of the Poincar\'e normal form. We define
\begin{align}
\label{zfromwn}
z = w + \frac{h_{20}}{2} w^2 + h_{11} w \bar{w} + \frac{h_{02}}{2} \bar{w}^2 + \frac{h_{30}}{6} w^3 +\frac{h_{21}}{2} w^2 \bar{w} +  \frac{h_{12}}{2} w \bar{w}^2 + \frac{h_{03}}{6} \bar{w}^3,
\end{align}
For Hopf bifurcations, it has been established that $h_{jk}$ can be chosen so that $\dot{w}$ depends, to lowest order, solely on $\lambda w$, and a single third order term:
 \begin{align}
  \label{wdotn}
 \dot{w} = \lambda w +c_1 w^2 \bar{w} + \mathcal{O}(|w|^4).
  \end{align}
That is, all other second- and third-order terms vanish except for a single third order term $w^2\bar{w}$, with coefficient $c_1$ \cite{Bautin1949} given by
     \begin{align}
  \label{c1n}
c_1 = \frac{2\lambda + \bar{\lambda}}{2|\lambda |^2} g_{20}g_{11} + \frac{1}{\lambda} |g_{11}|^2 + \frac{1}{2(2\lambda - \bar{\lambda})} |g_{02}|^2 + \frac{1}{2}\tilde{g}_{21}
  \end{align}
 The values $h_{jk}$ that achieve this cancellation in (\ref{wdotn}) are related to the $g_{jk}$ in the same way as in the $n=2$ case, except $\tilde{g}_{21}$ in the $n$-dimensional case has extra terms shown in (\ref{zdotnG21}). 
 
 Now we are setup to present the $n$-dimensional version of the theorem.
 
{\bf theorem}: The Hopf Mean Deviation Theorem for $n$-Dimensional Systems 

Consider a smooth $n$-dimensional system $\dot{x} = f(x; \alpha)$ undergoing a codimension one Hopf bifurcation at $\alpha^*$, inducing an limit cycle solution $x(t)$ for $\alpha \in \tilde{I}_{\alpha^*}$ emerging from an equilibrium $x_0(\alpha)$. The oscillatory mean $\langle x \rangle_{\alpha} = \frac{1}{T} \int_0^T x(t) dt$ deviates from the equilibrium according to
  \begin{align}
 \label{mainResultTn}
 \langle x \rangle_{\alpha} -  x_0(\alpha)  =  \begin{cases}  K(\alpha) \mu(\alpha)  + \mathcal{O}(\mu(\alpha)^2), &  \alpha \in \tilde{I}_{\alpha^*} \\ 0,  &  \alpha \in I_{\alpha^*}  \end{cases}.
\end{align}
The vector quantity $K$ in (\ref{mainResultTn}) is given by
\begin{align}
   \label{Kfinaln}
  K =   - Re \Big( 2 \frac{g_{11}}{ \bar{\lambda} Re(c_{1}) }\Big) Re(q)  + Im \Big( 2 \frac{g_{11}}{ \bar{\lambda} Re(c_{1}) } \Big)  Im(q)  - \frac{1}{Re(c_1)} \eta_{11} ,
  \end{align}
  where $g_{11}$, and $c_{1}$ are defined in (\ref{gcoef}), (\ref{c1n}), respectively.

{\bf Proof}
The circular limit cycle solution of (\ref{wdotn}) is $w(t) = r_w e^{i \omega_w t}$ with $r_w = \sqrt{-\mu/Re(c_1)}$ and $\omega_w = \omega - \frac{Im(c_{1})}{Re(c_{1})}\mu$. We use the fact that  $\int_0^T w^j(t)\bar{w}^k(t) dt  = 0$ when $j \neq k$, and $\int_0^T w(t)\bar{w}(t) dt  = r_w^2 = -\mu/Re(c_1)$, when $j = k = 1$, to compute, to lowest order, the mean of $z(t)$
    \begin{align}
  \label{zmeann}
\langle z \rangle = \frac{1}{T}\int_0^T z(t) dt =   -\frac{g_{11}}{\bar{\lambda}Re(c_1)} \mu + \mathcal{O}(\mu^2),  
 \end{align}
 and the mean of $V(z, \bar{z})$
 \begin{align}
\label{Vmean}
\langle V(z,\bar{z} ) \rangle  = -\eta_{11} \frac{1}{Re(c_1)}\mu  + \mathcal{O}(\mu^2).
  \end{align}
  The combined results of (\ref{zmeann}-\ref{Vmean}) enable the averaging of $x(t)$ on the invariant unstable manifold via (\ref{Wlock}) and (\ref{yWc}):
  \begin{align}
  \langle x\rangle  = \langle z \rangle q + \langle \bar{z} \rangle \bar{q} +  \langle V(z,\bar{z} )  \rangle + \ldots .
  \end{align} 
 The desired cycle mean of $x(t)$ is 
    \begin{align}
   \label{kProof}
\frac{1}{T}\int_0^T x(t) dt - x_0(\alpha) =  K \mu + \mathcal{O}(\mu^2),
  \end{align}
  where the vector quantity $K$ is given by
    \begin{align}
   \label{Kfinaln}
K =   - Re \Big( 2 \frac{g_{11}}{ \bar{\lambda} Re(c_{1}) }\Big) Re(q)  + Im \Big( 2 \frac{g_{11}}{ \bar{\lambda} Re(c_{1}) } \Big)  Im(q)  -\eta_{11}\frac{1}{Re(c_1)} .
  \end{align}

This $n$-dimensional version of the theorem contains the same terms involving $g_{11}$ as the $n=2$ version, but includes an additional term involving the vector $\eta_{11}$ that accounts for the distortion of the unstable manifold $W^u_{loc}(x_0)$ from $T^u$. 
 
$\square$
 
\subsection{Mean deviation in a $3$-dimensional example: a feedback control system}
\label{FCS}
Consider the $n=3$ dimension feedback control system with $x = (x_1,x_2,x_3)$:
\begin{align}
   \label{fbcs}
   \begin{split}
   \dot{x}_1 & = x_2, \\
   \dot{x}_2 & = x_3, \\
   \dot{x}_3 & = -\alpha x_3 - \beta x_2 +x_1(x_1-1),
   \end{split}
  \end{align}
  where $\alpha$ and $\beta$ are positive parameters. The equilibrium is the origin $x_0 = (0,0,0)$ for all parameter values. The Jacobian at the equilibrium has the characteristic equation
  \begin{align}
  \begin{split}
  \lambda^3+ \alpha \lambda^2 + \beta \lambda + 1 = 0.
  \end{split}
  \end{align}
 If we take $\beta$ as a fixed parameter and $\alpha$ as the parameter of bifurcation, we find that the equilibrium is linearly stable when $\alpha > \beta^{-1}$, and the bifurcation point is $\alpha^* = \beta^{-1}$, wherein $\lambda_j = \pm i \sqrt{\beta}$ for $j = 1,2$, respectively, and $\lambda_3 = -\beta^{-1}$. We have chosen $\beta = 1$ for this example. After finding the right $q$ and left $p$ eigenvectors associated with $\lambda_{1}$ and $\bar{\lambda}_1$, respectively, we compute the second-order tensors to be
 \begin{align}
 B(q,q) = \begin{bmatrix} 0 \\ 0 \\ 2q_1^2 \end{bmatrix}, \quad B(q,\bar{q}) = \begin{bmatrix} 0 \\ 0 \\ 2|q_1|^2 \end{bmatrix}.
 \end{align} 
The third order tensor is zero $C(x,y,z) = 0$ in this particular example. It is then straightforward to compute the $g_{jk}$ defined in (\ref{gcoef}), the $h_{jk}$ in (\ref{hquadg}-\ref{hcubic}), and the $H_{jk}$ in (\ref{HjkDefn}). We then use the $H_{jk}$ and the Jacobian to compute the $\eta_{jk}$ vectors that make up the lowest order quadratic components of $y = V(z,\bar{z})$ in (\ref{yWc}) by solving (\ref{Heqns}). We use the $g_{jk}$ and compute $\tilde{g}_{21}$ in (\ref{zdotnG21}), and then use $\tilde{g}_{21}$ and other coefficients to compute $c_1$ (\ref{c1n}). Combining $c_1$, $g_{11}$, and $\eta_{11}$, we can map the circular limit cycle solution $w(t) = r_w e^{i \omega_w t}$ into the $x(t)$ phase space via (\ref{zfromwn}), (\ref{Wlock}), and  (\ref{yWc}).

\begin{figure*}
  \centering
  \includegraphics[width=15cm]{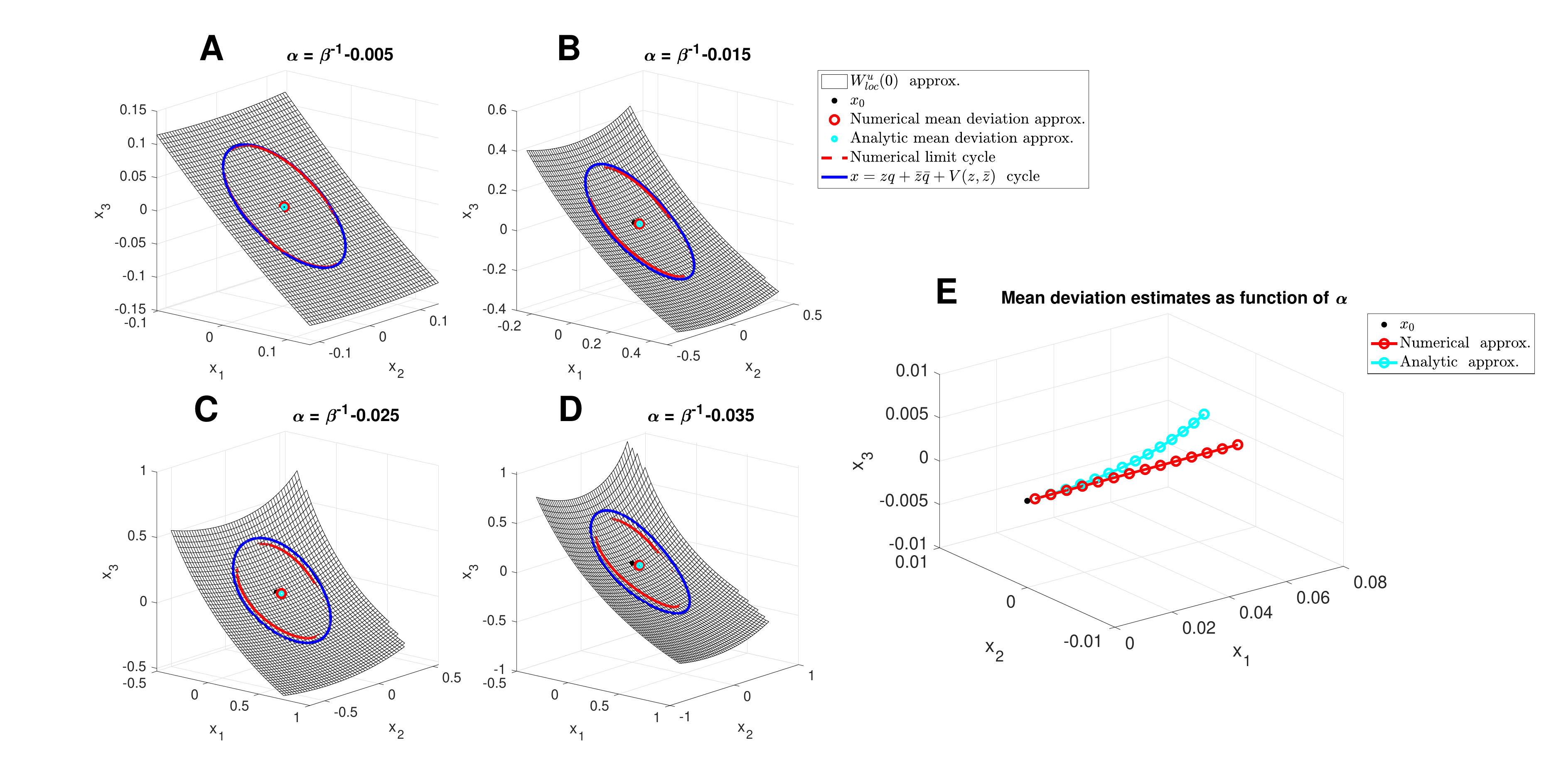}
  \caption{Mean deviation in an $n=3$-dimensional feedback control model. (A-D): Each panel shows the 3-dimensional phase-space, including the quadratic approximation of the $W^u_{loc}(x_0)$ surface, the numerical and analytic mean deviations, and numerical and analytic limit cycle orbits. A-D show these objects over four $\alpha$-values, each progressively further from the bifurcation point. Note the difference in scales between panels. (E): A close up of the numerically computed and analytically computed mean deviation $K(\alpha) \mu(\alpha)$ over a larger sampling of $\alpha$-values past the bifurcation point. Naturally, nearest to the bifurcation, when $\mu(\alpha)$ is near zero, the analytical approximation is tangent to the numerically computed mean deviation.}
  \label{Fig5}
\end{figure*}

Figure \ref{Fig5} A shows the phase space of the 3-dimensional system for an $\alpha$-value just past the bifurcation point where a small-amplitude stable oscillation exists. The quadratic approximation of the unstable manifold surface $W^u_{loc}(0)$ (using \ref{yWc}) passes through equilibrium at the origin. The numerical and analytic cycle mean deviations are also depicted very close to the origin. The numerical solution is well-approximated by the projection of the $w(t) = r_w e^{i \omega_w t}$ solution into the phase space. Figures \ref{Fig5} B-D show the numerically simulated orbit at increasingly large amplitudes (note the distinct axis scales for each A-D) for $\alpha$-values progressively further from the bifurcation point. Here, $W^u_{loc}(0)$ appears progressively more distorted as the quadratic terms of (\ref{yWc}) increasingly characterize the curved shape of the manifold. Additionally, the analytic estimate of the limit cycle becomes progressively less accurate at predicting the numerically simulated limit cycle solution. 

Finally, the numerical and analytic cycle means, owing to their $\mathcal{O}(\mu)$-magnitude dependence, increase their deviation from the origin as $\alpha$ moves further past the bifurcation. Figure \ref{Fig5}E shows a close-up of the numerically computed and analytic cycle mean deviations (\ref{kProof}-\ref{Kfinaln}) over a finer sampling of $\alpha$-points past the bifurcation. For small $\mu(\alpha)$ near the bifurcation we observe the $\mathcal{O}(\mu)$ analytic approximation is tangent to the numerical cycle mean deviation, but increasing $\alpha$ further yields greater divergence between the two estimates of the mean deviation, consistent with what was observed in the 2D examples.

\section{Discussion}
\label{discussionSection}


The study of how quantities derived from dynamical systems are affected by parameters---i.e., $\alpha$---has broad applications in science and engineering. Our mean value theorem contributes additional feature that can be predicted when a parameterized dynamical system is near a codimension one Hopf bifurcation. Namely, we have established that mean of a dynamical system as a function of $\alpha$ can deviate from the equilibrium $x_0$ when $\alpha$ crosses the bifurcation threshold $\alpha^*$, in which  $\langle x \rangle_{\alpha} =  x_0(\alpha) +  K(\alpha) \mu(\alpha) + \mathcal{O}(\mu(\alpha)^2)$. We have documented in several examples that the $K(\alpha) \mu(\alpha)$ term can dramatically alter the trajectory of the mean $\langle x \rangle_{\alpha}$ away from the equilibrium. 

Our mean value deviation result utilized the $g_{jk}$ coefficients of the dynamical system $f$, derived from the first- through third-order tensors of $f$ expanded about the equilibrium. The $g_{11}$ coefficient and the vector $\eta_{11}$ (in $n$-dimensional systems) are of particular importance to the mean deviation. This order-$\mathcal{O}(\mu)$ scaling result complements the longstanding knowledge about Hopf bifurcations, that the oscillation amplitude $r$ of the Hopf-induced limit cycle has an order-$r \sim \mathcal{O}( |\mu|^{1/2})$ scaling, as is detailed in \cite{kuznetsov98}.

There is another distinct approach to obtain the standard $r \sim \mathcal{O}( |\mu|^{1/2})$-amplitude orbit characterization of Hopf-induced cycles termed the `method of averaging' \cite{marsden2012hopf}. In this approach, the limit cycle is expressed as $z(t) = re^{i \theta}$. The oscillation amplitude $r$ is expressed as a periodic function $r(t,\theta,\alpha)$, in which $r$ depends on $\theta$ in a $2\pi$-periodic manner that can in effect shift the centroid of $z(t)$ from the origin. In this prior treatment, a $2\pi$-periodic correction term $u \sim \mathcal{O}(\mu(\alpha)^{1})$ is applied ($r \to r+u$) so that orbit is re-centered to the origin and the corrected average is zero $\int_0^{2\pi}r d\theta = 0$ (see page 155 in Section 4C of \cite{marsden2012hopf} for details). The correction term $u$ was devised by averaging over the relevant lowest-order polynomial expansions of the dynamics, which is similar to our approach in this article. However, while the existence of such a correction was established in this previous work, it was only used to recenter the orbit. The significance of this mean correction was never interrogated further. In particular, there was not an in-depth analysis of the contribution of individual low-order tensor components of $f$ to this mean correction $u$.  In contrast, the approach we have taken here demonstrates the significance of $g_{11}$ and $\eta_{11}$ for determining the mean deviation theorems (\ref{mainResultT}) and (\ref{mainResultTn}). 

By interpreting $\alpha$ as an input to the system, and the mean value of $x$ (\ref{xmean0}) as the output, we have identified that this mean deviation can produce what is termed oscillation-induced gain modulation---OIGM. In the predator-prey model (Section \ref{PPexample}) we have established that increasing $\alpha$ past the bifurcation threshold resulted in drop in the mean level of predators relative to the location of the predator-prey equilibrium (Figs \ref{Fig1}-\ref{Fig2}). 

Furthermore, the Brusselator model (Section \ref{BrusModel}) demonstrated that a wide variety of OIGM behaviors can be exhibited through changing a secondary parameter.  By changing a secondary parameter, we showed a Hopf bifurcation could alter the slope of the a mean quantity, as a function of $\alpha$, to shift upwards, stay the same, or shift downwards. Moreover, all of these OIGM behaviors were accurately predicted by our mean deviation theorem. 

We have also documented a novel OIGM phenomenon in the Wilson-Cowan model of coupled excitatory and inhibitory (e-i) neuronal populations \cite{wilson1972} (Section \ref{WCModel}). When input current $I$ were below the Hopf point and a stable equilibrium existed, there was steeper input gain in the mean excitatory activity, whereas the mean inhibitory activity showed a less-steep gain. However, when $I$ surpassed the bifurcation threshold, the relative gain between excitatory and inhibitory activity reversed: the inhibitory population exhibited higher gain than the excitatory population. This role-reversal induced by Hopf-induced oscillations could have functional consequences for excitatory-inhibitory balance within large-scale networks \cite{hennequin18}. Furthermore, any postsynaptic targets that are sensitive to the mean output of these e-i networks will necessarily be affected by the presence of these Hopf-induced oscillations.

In this article, we have stuck to relatively simple model examples found in standard textbooks. This was done for pedagogical reasons, but also to establish the generality of this result. We anticipate that more elaborate models that are used in contemporary applications will find this result useful as well.



\bibliographystyle{siamplain}
\bibliography{references}

\end{document}